\documentclass[10pt]{article}
\usepackage{amsmath,amsthm,amsfonts,amssymb}
\usepackage{graphicx}
\usepackage{hyperref}
\usepackage[usenames]{color}
\usepackage{amsfonts,amsmath, amsthm, amscd, amssymb}
\usepackage{diagrams}

\newtheorem{theorem}{Theorem}[section]
\newtheorem{lemma}[theorem]{Lemma}
\theoremstyle{definition}

\newtheorem{corollary}[theorem]{Corollary}

\newtheorem{definition}[theorem]{Definition}

\newtheorem{remark}[theorem]{Remark}
\newtheorem{proposition}[theorem]{Proposition}

\newcounter{comcount}

\newcommand{\factor}[2]{{\raise0.7ex\hbox{$#1$} \!\mathord{\left/
 {\vphantom {#1 {#2}}}\right.\kern-\nulldelimiterspace}
\!\lower0.7ex\hbox{${#2}$}}}

\def\MZ{{\mathbb{Z}}}

\title{Krull dimension of solvable groups}

\author{A. Myasnikov, N. Romanovskiy}

\date{ July 22, 2008}

\begin{document}
\maketitle

\begin{abstract}
In this paper we prove that free solvable groups have finite Krull dimension. In fact, this is true for much wider class of solvable groups, termed  rigid groups. Along the way we study the algebraic structure of the limit solvable groups (fully residually free solvable groups).
\end{abstract}

\tableofcontents

\section{Introduction}
\label{se:intro}

In this paper we study the Krull dimension of solvable groups. We show that the Krull dimension is finite for a wide class of solvable groups, including polycyclic groups, free solvable groups of finite rank, and iterated wreath products of free abelian groups of finite rank. Free solvable groups and iterated wreath products of torsion-free abelian groups share some common rigidity properties. To  study these and similar groups  we introduce a new class of solvable groups, termed rigid groups. The Diophantine geometry in rigid groups is tractable and the Krull dimension is finite. Furthermore, our approach gives numerical estimates on the Krull dimension of a rigid group $G$ in terms of some natural  algebraic invariants of $G$. Along the way we lay down foundations of dimension theory in groups, tying down Krull, Zariski, and projective  dimensions altogether.  As a bi-product we  describe algebraic structure of limits of rigid groups $G$ in Gromov's topology, i.e., fully residually $G$-groups.

Dimension is one of the fundamental notions in geometry and topology, while its algebraic counterpart, the Krull dimension,  plays an equally fundamental roll in commutative algebra and algebraic geometry.  The classical theory of topological dimension goes back to the works of Poincare and Brouwer.  Now there are several different notions  of dimension of a topological space,  among which  the {\em Menger dimension} and the {\em Lebesgue covering dimension} are the most familiar ones (see  \cite{HW,Pears} for details). In algebraic geometry the typical topological spaces are Noetherian, in which case the notion of Zariski dimension  arises naturally \cite{Hartshorne,Shafarevich,Eizenbud} and  gives   the dimension of irreducible  varieties in affine spaces. In general, all these dimensions are different, though in many cases they coincide.  Here we approach  groups via algebraic or Diophantine geometry over groups (see \cite{BMR,MR2} for details), so our primary concern is the Zariski dimension.  It turns out that for a group $G$ the Zariski topology   on $G^n$ (for each $n$) is Noetherian if and only if the group $G$ is equationally Noetherian, i.e., every system of equations in finitely many variables and coefficients from $G$  is equivalent in $G$ to a finite subsystem of itself \cite{BMR}. This result gives a convenient tool to approach Zariski topology in groups in a pure algebraic way. During the last several years the class of equationally Noetherian groups was intensively studied, many groups appeared to be in this class: linear groups \cite{Bryant,Guba,BMR},  hyperbolic groups \cite{Sela7}, some interesting classes of solvable groups (including free solvable groups) \cite{GR}. This opens the way to study the Zariski dimension in these groups.

 On the algebraic side it was known already in the beginning of  the twenties century that the Zariski dimension of a variety $Y$ in an affine space over  a field $k$  equals to the transcendence degree of the field of the rational functions on $Y$ over $k$.  In 1940's Krull introduced a new type of a dimension (now known as the Krull dimension) in an arbitrary  ring $R$  as the supremum of the heights of prime ideals of $R$, and showed that for the coordinate algebras $A(Y)$ of irreducible affine varieties $Y$ the new dimension is equal to the Zariski one for $Y$. Nowadays the Krull dimension is one of the fundamental notions in ring theory, we refer to \cite{Eizenbud} for a comprehensive treatment of the Krull dimension. In this paper by analogy with commutative algebra we introduce a notion of a prime ideal in a group, which generalizes the one from \cite{BMR}, and define a Krull dimension, which is more easier to deal with in the algebraic context than the Zariski  dimension. Bringing algebra and geometry together  we prove that the Zariski dimension of an algebraic irreducible set $Y$ in $G^n$ is equal to the Krull dimension of the coordinate group $\Gamma(Y)$. We use the Krull dimension in our  treatment of dimension in solvable groups, but in a slightly different, though equivalent form, that is especially suitable in the context of combinatorial group theory.  This variation corresponds to the dimension of algebraic sets defined by equations without coefficients. In algebraic terms it can be described as the projective dimension (quotient height) of a group $G$, which is the supremum of the lengths of chains of proper epimorphism of groups $G = G_1 \to G_2 \to \ldots \to G_{k+1}$, where $G_i$ are universally equivalent to $G$.

Krull dimension is one of the crucial numerical invariants of groups. Computing the dimension in a given group is a hard task, so concrete results are scarce here.  To mention a few, observe, first, that for a free abelian group $A$ the Zariski dimension of $A^n$ is equal to $n$. Polycyclic groups provide another example of groups with finite Zariski or Krull dimension.
Recall, that a polycyclic group $G$ is linear, hence the Zariski  topology on $G^n$ is Noetherian for each $n$ \cite{BMR}. If $G$ is torsion-free then  the Krull dimension of $G$ is finite and bounded from above by the Hirsch length of $G$ (see Proposition \ref{pr:polycyc}), which is the number of infinite cyclic factors in a polycyclic series of $G$ (it does not depend on the series).  However, the precise values of the Krull dimension for torsion-free polycyclic groups are unknown. Finitely generated metabelian groups are also equationally Noetherian \cite{BMR}, so the Zariski dimension is well defined in such groups.  In a pioneering  paper \cite{Rem-Tim} Remeslennikov and Timoshenko calculated the projective dimension  of some metabelian  groups, including free metabelian groups of finite rank and wreath products of free abelian groups of finite rank. This rests on  a thorough  study of groups that satisfy all universal sentences of the theory of  free metabelian groups (called $u$-groups) \cite{Chapuis,Rem-Shtor2}.  A detailed description of algebraic sets in $u$-groups is contained  in a series of papers \cite{Rem-Shtor1,Rem-Rom1,Rem-Rom2}.

In this paper we introduce a new class of solvable groups, termed rigid groups, where the Krull dimension (as well as the projective dimension) is finite. The class of rigid groups contains free solvable groups and iterated wreath  products of torsion-free abelian groups. Furthermore, we give numerical estimates of the dimensions involved. Along the way we describe the algebraic structure of groups universally equivalent  to a given rigid group. In particular, we show that these groups are again rigid.

The unification theorem mentioned above implies that the coordinate groups of irreducible varieties in rigid groups are also rigid, as well as the limits of rigid groups in the Gromov-Hausdorff topology \cite{Gromov}, see also \cite{Grigorchuk}, and more details in \cite{Cham,CG}. Notice, that the limit groups of free groups (known also as fully residually free groups, or freely discriminated groups, or $\omega$-residually free groups) played a crucial part in a recent solution of the Tarski problems \cite{KM4,Sela6}, now they are one of the most actively studied classes of groups in geometric and combinatorial group theory (see, for example, \cite{Bridson1,Bridson2,BHMS,Groves1,Groves2,G1,CZ} and the bibliography there). The place of rigid groups in the theory of solvable groups is quite similar to the one of the limit groups in the class of all groups. It seems plausible  that our approach to rigid groups could be used to develop a theory of limit groups in the category of solvable groups.

\section{Irreducible algebraic sets and their coordinate groups}
\label{se:coordinate-groups}

In this section we discuss some basic notions on irreducible algebraic sets in groups and refer to papers  \cite{BMR,MR2} for details.

Let $G$ be a group  and  $F(X)$ a free
group with basis $X = \{x_1, x_2, \ldots  x_n\}$. By  $G[X] = G \ast
F(X)$ we denote  a free product of $G$ and $F(X)$. If $s \in G[X]$ then
the formal expression $s = 1$ is called an {\em equation}
over $G$, more generally, for   $S \subset G[X]$
 the expression $S = 1$ is called a {\em system of equations}
over $G$. In an equation $s = 1$ (or a system $S =1$) the element $s$, as an element of the free product $G[X]$, can be written as a product of some
elements from $X \cup X^{-1}$, termed {\em variables},
and some elements from $G$, termed   {\em constants}. To indicate  variables in a system $S = 1$ we
sometimes write  $S(X) = 1$. If $s \in G[X]$ and $g_1, \ldots, g_n \in G$ then by $s(g_1, \ldots, g_n)$ we denote the element of $G$ obtained from $s$ by replacing each $x_i$ by $g_i$, $i = 1, \ldots,n$.
A {\em solution} of the system $ S(X) = 1$ in $G$ is a
tuple of elements $(g_1, \ldots, g_n) \in G^n$ such that  $s(g_1, \ldots, g_n) = 1$ for every $s \in S$.

 The  set
$V_G(S)$ (or simply $V(S)$)   of all solutions of the system $S(X)=1$ is called an {\em
algebraic subset} in  $G^n$. Two systems $S(X) = 1$ and $T(X) = 1$ are {\em equivalent} in $G$ if $V_G(S) = V_G(T)$. There exists the largest (relative to inclusion) subset  $R(S) \subseteq G[X]$ such that the system $R(S) = 1$  is equivalent to $S(X) = 1$:
 $$R(S) = \{f \in G[X] \ | \ f(g_1, \ldots, g_n) = 1 \   \forall
(g_1, \ldots, g_n) \in V_G(S) \},$$
which is termed  the {\em radical } of $S$. Clearly, $R(S)$  is a normal subgroup of $G[X]$.
Notice that if
$V_G(S) = \emptyset$, then $R(S) = G[X]$. Sometimes the following notation is in use. If $Y$ is an algebraic set in $G^n$ then the set
$$I(Y) = \{f \in G[X] \mid \ f(a) = 1 \   \forall a \in Y \}$$
 is a normal subgroup of $G[X]$, moreover, if $Y = V_G(S)$ then $I(Y) = R(S)$.

The quotient group $\Gamma(Y) = G[X]/I(Y)$
 is the {\em coordinate group} of an
algebraic set $Y \subseteq G^n$.  If $Y = V(S)$  then, sometimes, we refer to the coordinate group $\Gamma(Y)$ as to $G_{R(S)}=G[X]/R(S)$.  It is not hard to show that if $Y$ is a non-empty algebraic set  then $G \cap I(Y) = 1$, so the group $G$ canonically  embeds into the coordinate group $\Gamma(Y)$.
Notice, that every solution $(g_1, \ldots, g_n)$ of $S(X) = 1$ in $G$
can be described as a homomorphism $\phi: G_{R(S)} \rightarrow G$ such that $\phi(x_i) = g_i, i = 1, \ldots,n, $ and $\phi(g) = g$ for $g \in G$. The converse is also true.

It is convenient to study coordinate groups in the framework of  $G$-groups. A {\em $G$-group} is an arbitrary group with a distinguished subgroup isomorphic to $G$. Most of the basic group-theoretic notions make sense for  $G$-groups. For example, if $H$ and $K$
are $G$-groups then a homomorphism $\phi: H \rightarrow K$ is a
{\em $G$-homomorphism} if $g^\phi = g$ for every $g \in G$; {\em $G$-subgroups} of $H$ are subgroups containing  the subgroup $G$, etc. It is easy to see that the free product $G[X] = G \ast
F(X)$ is a free object in the category of $G$-groups (with $G$-homomorphisms as the morphisms).   A $G$-group $H$ is termed {\em finitely generated
$G$-group} if there exists a finite subset $A \subset H$ such that
the set $G \cup A$ generates $H$. From now on when we refer to the coordinate group $\Gamma(Y)$ for $Y \neq \emptyset$, we assume that the distinguished subgroup is the the image of $G$ under the canonical  embedding $G \rightarrow \Gamma(Y)$ described above. In this terminology, the set of solutions of $S(X) = 1$ in $G$ can be described as the set of $G$-homomorphisms from $G/R(S)$ onto $G$. We refer to \cite{BMR} for a
general discussion on $G$-groups.

Recall that a group $G$ is {\em equationally Noetherian} (\cite{BMR,MR2,BMRom}) if for any finite set $X$ any
system of equations from $G[X]$ is equivalent over $G$ to some finite part of
itself.     Equationally   Noetherian groups enjoy  properties similar to Noetherian rings.
Among the most interesting examples of  equationally  Noetherian groups we would like to mention the following ones: linear groups over Noetherian rings (in particular, free groups,
polycyclic groups, finitely generated metabelian groups) \cite{Bryant,Guba,BMR}, torsion-free
hyperbolic groups \cite{Sela7}, and free solvable groups \cite{GR}. It has been shown in \cite{BMR} that a group $G$ is equationally Noetherian if and only if for every $n$  the Zariski  topology on $G^n$ is {\em  Noetherian}, i.e., every proper descending chain of closed sets in $G^n$ is
finite. Recall that the {\em  Zariski topology} on $G^n$ is defined by taking the algebraic sets as a pre-basis of closed sets.

It is easy to see that if the Zariski topology is  Noetherian then every closed  set $V$ in $G^n$ is a finite union of irreducible pair-wise distinct closed subsets (called {\it irreducible components} of $V$):
$$V = V_1 \cup \ldots V_m.$$
Under requirement $V_i \not \subseteq V_j$ for $j \neq j$ the components of $V$ are uniquely determined.
 Recall that a closed
subset in $G^n$ is {\it irreducible} if it is not a union of two
closed proper subsets.

 To formulate the principal  result on the coordinate groups of irreducible algebraic sets we need a few definitions.

 A $G$-group $H$  is  termed {\em $G$-discriminated} by $G$ if for any finite subset $K \subseteq H$ there exists a $G$-homomorphism $\phi:H \to G$ such that the restriction of $\phi$ on $K$ is injective. If the condition above holds only for singletons $K$ then the group $H$ is called {\em $G$-separated} by $G$. Sometimes, we omit $G$ and say that $H$ is  discriminated by $G$.

Let ${\cal L}$ be the group theory language, consisting  of a binary operation
$\cdot$ (multiplication), a unary operation $^{-1}$ (inversion),
and a constant symbol $1$ (the identity). By ${\cal L}_G$ we denote the extended language ${\cal L} \cup \{c_g \mid g \in G, g \neq 1\}$ where all non-trivial elements from $G$ are added to ${\cal L}$ as new constants.
 We say that two $G$-groups $H$ and $L$ are {\em $G$-universally equivalent} ({\em $G$-existentially equivalent})  if they satisfy the same universal (existential) first-order sentences in the language ${\cal L}_G$.
Similar definitions apply in the case of rings, or any other category of algebraic structures.

\begin{theorem}
\label{coord-groups} \cite{BMR,MR2}

 Let $G$ be an equationally Noetherian group. Then for a finitely
generated $G$-group $H$ the following conditions are
equivalent:
\begin{enumerate}
\item [1)] $H$ is $G$-universally equivalent to $G$;
\item [2)] $H$ is $G$-existentially equivalent to $G$;
\item [3)] $H$ $G$-embeds into an  ultrapower of $G$;
\item [4)] $H$ is $G$-discriminated by $G$;
\item [5)] $H$ is the  coordinate group of an irreducible algebraic
set over $G$ defined by  a system of equations with coefficients
in $G$.
\end{enumerate}

\end{theorem}

\section{Krull dimension in groups}
\label{se:Krull}

In this section we discuss the Krull dimension in groups.

The {\em dimension} (see, for example,  \cite{Hartshorne,Shafarevich,Eizenbud})  of a Noetherian topological space $T$ is defined as the supremum of the natural numbers $k$ such that there exists a chain
$$Z_0 \subset Z_1 \subset \ldots  \subset Z_k$$
 of distinct irreducible non-empty closed subsets of $T$. This dimension is also known as   the  {\em Zariski dimension} (see, for example, \cite{Z,Hrush}) or the {\em combinatorial dimension}  \cite{Matsumura}.  In this paper we call it  Zariski  dimension.

If $G$ is an equationally Noetherian group then an  algebraic subset $Y$ of $G^n$ is a Noetherian topological space in the Zariski topology induced from $G^n$, so its Zariski dimension is well-defined. We say that an equationally Noetherian group $G$ has a finite Zariski dimension if
for every $n$ the Zariski dimension of  $G^n$ is finite.

One of the classical results in algebraic geometry relates the Zariski dimension of an affine algebraic set $Y$ with the Krull dimension of the coordinate algebra of $Y$. Namely, let $k$ be an algebraically closed field and   $k[X]$ the $k$-algebra of polynomials in variables $X = \{x_1, x_2, \ldots  x_n\}$ and coefficients in $k$. A subset   $S \subseteq k[X]$ defines an affine algebraic set
$$V(S) = \{a \in k^n \mid f(a) = 0 \ \text{for all } a \in V(S)\}.$$
  By Hilbert Basis theorem the Zariski topology on the  affine space $k^n$ (with algebraic sets as closed subsets) is Noetherian, so the Zariski dimension of $Y = V(S)$ is well-defined (and finite). The  {\em coordinate  algebra}  $A(Y)$ of an algebraic set $Y \subseteq k^n$ is the quotient $k[X]/I(Y)$, where
  $$
  I(Y) = \{f \in k[X] \mid f(a) = 0 \ \text{for  all } \ a \in Y\}.
  $$
is an ideal in $k[X]$ (denoted also by $R(S)$).    The crucial result tying the geometric and algebraic view-points on dimension states that  the Zariski  dimension of $Y$ is equal to the Krull dimension of $A(Y)$. Here,  the {\em Krull dimension} of a commutative unitary ring $A$ is the supremum of  all natural numbers $k$ such that exists a chain
$$p_0 \subset p_1 \subset \ldots  \subset p_k$$
of distinct prime ideals in $A$. We refer to  \cite{Eizenbud} for  details on Krull dimension.

It turns out that a similar result holds in the group case as well. The corresponding notion of a prime ideal in a $G$-group was introduced in \cite{BMR}. Here we   use another,  more general definition, which coincides with the original one for all groups discussed in \cite{BMR}. To explain we digress first into commutative algebra.

 It is easy to see that (in the notation above) the ideal  $R(S)$ has a special property, it is the intersection of the kernels of all $k$-homomorphisms $\phi:k[X] \to k$ such that $\phi(S) = 0$.   The ideal $p$ in $k[X]$ is called a {\em radical ideal} or {\em closed} ideal if $p = R(p)$. By Hilbert Nullstellensatz
$$R(p) = \{f \in k[X] \mid f^n \in p \ \text{for  some  positive  integer }  n \},$$
so every prime ideal $p$ of $k[X]$ is radical and the algebraic set $V(p)$ is irreducible. The converse is also true, namely,  if $V(S)$ is an irreducible algebraic set then the ideal $R(S)$ is prime.  It follows that in the case when $k$ is algebraically closed an ideal $p$ of $k[X]$ is prime if and only if the $k$-algebra $k[X]/p$ is $k$-discriminated  by $k$.  Indeed, assume $p$ is prime and $y_1, \ldots, y_m$ are elements of $k[X] \smallsetminus p$, then their product $u = y_1 \ldots y_m$ is also not in $p$. Since $p = R(p)$ there exists $\phi:k[X] \to k$ such that $\phi(p) = 0$ and $\phi(u) \neq 0$, so $\phi$ factors through $k[X]/p$ and $\phi(y_i) \neq 0$ for each $i$. To prove the converse observe, that if $k[X]/p$ is $k$-discriminated by $k$ and $u,v \in k[X]/p$ are non-zero elements  with $uv = 0$ then there is a $k$-homomorphism $\phi: k[X]/p \to k$ such that $\phi(u) \neq 0, \phi(v) \neq 0$. But then $\phi(u)\phi(v) = 0$ - contradiction with $k$ being a field.

In general (when $k$ is not necessary algebraically closed),  $k[X]/p$  is discriminated by $k$ if and only if the ideal $p$ is prime and radical.  This characterization of prime and radical ideals holds for an arbitrary finitely generated commutative unitary $k$-algebra $A$. In this case by a radical ideal we understand an ideal $p$ of $A$ which  is  the intersection of the kernels of all $k$-homomorphisms $\phi:A \to k$ such that $\phi(p) = 0$.

Now it is clear how one can introduce the notion of prime radical ideals in groups.  We say that a normal subgroup $N$ of a $G$-group $H$ is an {\em ideal} in $H$ if  $N \cap G = 1$. An ideal $N$ is  {\em prime radical}, or just {\em prime},   if  $H/N$, as a $G$-group, is  $G$-discriminated by $G$.

\begin{lemma} \label{le:prime-ideal}
Let $G$ be an equationally Noetherian group, $X = \{x_1, x_2, \ldots  x_n\}$ a finite set of variables, and $p$ an ideal in $G[X]$. Then the algebraic set $V_G(p)$ is irreducible if and only if $R(p)$ is a prime ideal in $G[X]$.
\end{lemma}
\begin{proof}
Follows immediately from Theorem \ref{coord-groups}.
\end{proof}

Now, analogously to the classical case we define the {\em Krull dimension} $Kdim_G(H)$ of a $G$-group $H$ as the supremum of  all natural numbers $k$ such that exists a chain
$$p_0 \subset p_1 \subset \ldots  \subset p_k$$
of distinct prime ideals in $H$.

\begin{proposition} \label{prop:Krull}
Let $G$ be an equationally Noetherian group. If $Y$ is an algebraic subset of $G^n$ then the Zariski dimension of $Y$ is equal to the Krull dimension $Kdim_G(H)$ of its coordinate group $\Gamma(Y)$.
\end{proposition}
  \begin{proof}
  Straightforward from Lemma \ref{le:prime-ideal}.
  \end{proof}

 It is convenient to introduce another  of dimension. Let $H$ be a $G$-group discriminated by $G$. The {\em projective dimension} (or the {\em quotient height}) $pd_G(H)$ of $G$ is the supremum of  all natural numbers $k$ such that exists a chain
$$H =H_0 \to H_1  \to \ldots  \to H_k = G$$
 of proper $G$-epimorphisms of $G$-groups discriminated by $G$. If the group $H$ is not $G$-discriminated by $G$, then the projective dimension $pd_G(H)$ is defined as the supremum of projective dimensions of quotients of $H$ which are $G$-discriminated by $G$.

\begin{proposition} \label{prop:resolution}
Let $G$ be a  group and $H$ a finitely generated $G$-group $G$-discriminated by $G$. Then the Krull dimension of $H$ is equal to the projective dimension of $H$.
\end{proposition}
  \begin{proof}
  Straightforward from definitions.
  \end{proof}

Now we introduce one more variation of dimension  which is  suitable for combinatorial group theory. This variation corresponds to the dimensions of algebraic sets defined by equations without coefficients. To explain we look again at the definitions  in the case of $G$-groups:  if an ideal $N$ of a $G$-group $H$ is  prime then   $H/N$  is  $G$-discriminated by $G$.  Notice, that in this case $G$ is obviously  $G$-discriminated by $H/N$ (by the canonical embedding $G \to H/N$), so $G$ and $H/N$ are mutually $G$-discriminating each other. In terms of logic the former condition ($G$ $G$-discriminates $H/N$)  implies that $H/N$ satisfies all the universal sentences of group theory language $\mathcal{L}_G$  that hold in $G$, the so-called the {\em universal theory} $Th_{\forall,G}(G)$; while the latter one implies that $G$ and $H/N$ satisfy exactly the same universal sentences, i.e., $Th_{\forall,G}(G) = Th_{\forall,G}(H/N)$. If the group $G$ is equationally Noetherian and the group $H$ is finitely generated over $G$ then the converse is also true, so if $Th_{\forall,G}(G) = Th_{\forall,G}(H/N)$ then  $G$ and $H/N$ $G$-discriminate each other.  This shows that in general there are two possible definitions of a prime ideal, one - in terms of discrimination, and the other - in terms of universal equivalence.  In the coefficient-free case  there is an extra possibility, one may either require the mutual discrimination (universal equivalence) of groups at hands or does not.
  In fact, these discrepancy explains two different types of dimensions discussed in \cite{Rem-Tim}.

  Now we are in a position to define Krull dimension  of an abstract group $H$. We formulate it in terms of the {\em projective dimension} $pd(H)$, which is equal to the supremum of the lengths $k$ of chains of proper epimorphism of groups
   $$H = H_1 \to H_2 \to \ldots \to H_{k+1},$$
   where $H_i$ and $H$ mutually discriminate each other (without coefficients)  $i = 1, \ldots, k+1$. One can also consider a weaker property when in the definition above mutual discrimination is replaced by universal equivalence (which is implied by the mutual discrimination).  Since in this paper we consider only the case of equationally Noetherian finitely generated groups (or $G$-groups) then for us the both definitions are equivalent.

\begin{remark}
Let $H$ be a $G$-group $G$-discriminated by $G$. Then $pd_G(H) \leq pd(H)$.
\end{remark}

Now we show some known  examples of groups with finite Krull dimension.   For an abelian group $G$ the rank $rank(G)$ is defined as the cardinality of a maximal linearly independent set of elements in $G$.

\begin{lemma} \label{le:abelian-height}
Let $G$ be a free abelian group of rank $n$. Then $pd(G) = n$.
\end{lemma}
\begin{proof}
Notice, that any two torsion-free abelian groups are universally equivalent \cite{szmielew}. This implies that an ideal $N \leq G$ is prime if and only if $G/N$ is torsion-free (i.e., $N$ is a direct factor of $G$),  hence the result.
\end{proof}

Recall that a group $G$ is polycyclic if it has a finite subnormal series (called a {\em polycyclic series})
$$G = G_1 \unrhd G_2 \unrhd \ldots \unrhd G_k \unrhd G_{k+1} =  1$$
with cyclic quotients $G_i/G_{i+1}$, $i = 1, \ldots,k$.  The number $h(G)$ of  infinite cyclic quotients in a polycyclic series of $G$ does not depend on the series, it is called the Hirsch length of $G$.  Recall that polycyclic groups are linear, so the Zariski dimension is well-defined, and the Krull dimension, as well as the projective dimension, can be equivalently defined using the mutual discrimination or the universal equivalence.

\begin{proposition} \label{pr:polycyc}
Let $G$ be a torsion-free polycyclic group then $pd(G) \leq h(G)$.
\end{proposition}
\begin{proof}
In the proof we use some known facts on polycyclic groups, all of them can be found in  \cite{Kar-Mer}. Observe first, that if a group $G/N$ is universally equivalent to $G$ then $H = G/N$ is torsion-free and polycyclic.  In this event $H$ has a polycyclic series
$$H = H_1 \unrhd H_2 \unrhd \ldots \unrhd H_k \unrhd H_{k+1} =  1$$
with all the quotients $H_i/H_{i+1}$ infinite cyclic. Taking the full pre-images of the subgroups $H_i$ in $G$ under the canonical projection $G \to G/N$ gives a series
 $$G = G_1 \unrhd G_2 \unrhd \ldots \unrhd G_k \unrhd G_{k+1} =  N > 1$$
 in $G$ with the quotients $G_i/G_{i+1}$ infinite cyclic. The group $N$ is also polycyclic as a subgroup of $G$. Hence  combining the series above with the polycyclic series
 $$N = N_1 \unrhd N_2 \unrhd \ldots \unrhd N_m \unrhd N_{m+1} =  1$$
 of $N$ one gets a polycyclic series of $G$. Notice, that $N$ is finitely generated and torsion-free, so $h(N) \geq 1$. It follows that $h(G) \gvertneqq h(G/N)$. Hence the quotient length of $G$ cannot exceed its Hirsch length, as claimed.
\end{proof}

 If $M$ is a free metabelian group of finite rank $n$ then the projective dimension $pd(M)$   of $M$ is finite and its precise value is computed in   \cite{Rem-Tim}. In the same paper the projective dimension  of a discrete wreath product $W_{n,m}$ of two free abelian groups of ranks $m$ and $n$ was computed.    This shows that the Krull dimensions of $M$ and $W_{n,m}$ are finite.

 In the case of free or hyperbolic groups much less is known. If $F$ is a free group then the Zariski dimension of $F^1$  is equal to 2 \cite{Appel,Lorents,CR}. This gives the Krull dimension of the free product $F[X]$ (viewed as an $F$-group), where $X$ has cardinality one.

\section{Matrix splittings over a subgroup}

In this section we introduce some notation and prove several facts that are in use throughout the paper.
Some of the facts are known in folklore, but we add them to make the paper self contained.

Let $G$ be a group and $x,y \in G$. We denote by $x^y = y^{-1}xy$ the conjugation of $x$ by $y$,  and by $[x,y]$ the commutator $x^{-1}y^{-1}xy$ of $x$ and $y$.

Recall, that if $C$ is an abelian normal subgroup of $G$ then $G$ acts by conjugation on $C$ and this action turns $C$ into a right module over the group ring $\mathbb{Z}[G/C]$. Namely, for an element $u = m_1\bar{h}_1 + \ldots +m_k\bar{h}_k \in \MZ[G/C]$, where $m_i \in \MZ, \bar{h}_i = h_iC \in G/C, h_i \in G$, the right action of $u$ on $c \in C$ is defined (in the exponential notation) by
$$c^u = (c^{h_1})^{m_1} \ldots (c^{h_k})^{m_k}.$$
  Recall, that a right module $T$ over a ring $R$ has no $R$-torsion if for any $x \in T$ and $r\in R$ if $xr = 0$ then either $r= 0$ or $x = 0$.

\begin{lemma} \label{le:splitting}
Let $G$ be a group, $C$ an abelian  normal subgroup of $G$, $B = G/C$ and $--:G \to B$ the canonical epimorphism.  Then there exists a right $\MZ B$-module $D$ and a monomorphism $\lambda:G \to M(G)$ from $G$ into the group of matrices  $$M(G) = \begin{pmatrix} B & 0 \\ D & 1 \end{pmatrix},$$
such that for  $g \in G$  $$g\lambda = \begin{pmatrix} \overline{g} & 0 \\ d(g) & 1 \end{pmatrix},$$
 where $d(g) \in D$, and the following conditions hold:
    \begin{itemize}
    \item [1)] the module $D$ is generated by  $d(g), \ g \in G;$
    \item [2)] the map $d(g) \rightarrow \overline{g} -1$ extends to an epimorphism of  $\MZ B$-modules
     $$\sigma: D \to (B-1)\cdot \MZ B,$$
      where $(B-1)\cdot \MZ B$ is the fundamental ideal of the group ring  $\MZ B$,  viewed as $\MZ B$-module;
         \item [3)]  $\ker \sigma = d(C)$.
         \end{itemize}
 \end{lemma}
\begin{proof}
Let $G = F/N_1$ and $B = F/N_2$ be presentations of the groups $G$ and $B$, where $F= F(X)$ is a free group with basis $X = \{x_i \mid i \in I\}$. Denote the generators $x_iN_1$ of $G$ by $g_i$ and the generators $x_iN_2$ of $B$ by $b_i$, $i \in I$. Let  $T$ be a right free  $\MZ B$-module with basis  $\{ t_i \ \mid \ i \in I \}$  and
 $$\psi :T \rightarrow (B-1) \cdot \MZ B$$
 a module epimorphism which is defined by $t_i\psi = b_i - 1$, so  $(\sum t_iu_i) \psi = \sum (b_i-1)u_i.$

In \cite{Magnus} Magnus showed that the map
$$x_i \rightarrow \begin{pmatrix} b_i & 0 \\ t_i & 1 \end{pmatrix}, \ i \in I$$
defines a group homomorphism (Magnus homomorphism)
$$\varphi :F \rightarrow \begin{pmatrix} B & 0 \\ T & 1
\end{pmatrix}.$$
  with the kernel   $\ker \varphi = [N_2,N_2]$ (see also books \cite{Birman, Gupta}), so $\ker \varphi = [N_2,N_2] \leqslant N_1 \leqslant N_2$. Remeslennikov and Sokolov proved  in \cite{RS}  that
  $$(N_2)\varphi = \begin{pmatrix} 1 & 0 \\ U & 1 \end{pmatrix},$$
  where  $U= \ker \psi$.  It follows that
   $$(N_1) \varphi = \begin{pmatrix} 1 & 0 \\ U_1 & 1 \end{pmatrix},$$
   for some submodule  $U_1$ of  $U.$  This implies that $\varphi$ induces an  embedding $\lambda$ of $G$ into the following matrix group, where   $D=T/U_1$:
    $$\begin{pmatrix} B & 0 \\ D & 1 \end{pmatrix}.$$
  Clearly, the image of $C$ under this embedding is equal to
   $$\begin{pmatrix} 1 & 0 \\ U/U_1 & 1 \end{pmatrix},$$
     so it may be identified with the module $U/U_1.$

   The epimorphism  $\psi$ induces an epimorphism  $\sigma :D \rightarrow (B-1) \cdot ZB$ with the kernel $C.$  By construction one has $d(g_i) \sigma = b_i - 1= \overline{g_i}-1$ for the generators $g_i$ of $G$, which implies that $d(g) \sigma
=\overline{g}-1$ for all  $g \in G.$  Indeed, computing the image $(uv)\varphi$ for $u, v \in F$ one gets
  $d(uv) = d(u)\bar{v} + d(v)$, hence $d(uv)\sigma = (d(u)\bar{v} + d(v))\sigma = (\bar{u} - 1)\bar{v} +(\bar{v} - 1) = \overline{uv} - 1$, as claimed.
 This proves the lemma.
\end{proof}

In the situation above the group $M(G)$, together with the embedding  $\lambda:G \to M(G)$, is called a {\em matrix splitting} (or {\em splitting}) of $G$ over $C$. The following result shows that such   splittings are essentially unique.

\begin{lemma} \label{le:split-uniqueness}
Let $G$ be a group and $C$ an abelian  normal subgroup of $G$. Then the matrix splitting $M(G)$ of $G$ over $C$ is unique up to an isomorphism. Namely, if
$$\lambda_1: G \to M_1(G) = \begin{pmatrix} B & 0 \\ D_1 & 1 \end{pmatrix}$$
 is another splitting of $G$ over $C$ with
 $$g\lambda_1 = \begin{pmatrix} \overline{g} & 0 \\ d_1(g) & 1 \end{pmatrix},
\ g \in G,$$
then the map $d(g) \rightarrow d_1(g), \ g \in G,$ determines an isomorphism of modules $\tau: D \rightarrow D_1,$
which, together with the identity map  $B \rightarrow B$, gives rise to an isomorphism of the matrix splittings
$$\begin{pmatrix} B & 0 \\ D & 1 \end{pmatrix} \rightarrow \begin{pmatrix} B & 0 \\ D_1 & 1 \end{pmatrix}.$$
  \end{lemma}
\begin{proof}  We use notation from Lemma \ref{le:splitting}. Let $$\lambda_1: G \to M_1(G) = \begin{pmatrix} B & 0 \\ D_1 & 1 \end{pmatrix}$$
  be another splitting of $G$ over $C$. By definition  the module  $D_1$ is generated by the set   $\{ d_1(g_i) \ | \ i \in I \}.$ Hence the module homomorphism $\delta: T \to D_1$ that extends the map  $t_i \rightarrow d_1(g_i), \ i \in I,$ determines an epimorphism of groups
  $$\begin{pmatrix} B & 0 \\ T & 1 \end{pmatrix} \rightarrow \begin{pmatrix} B & 0 \\ D_1 & 1
\end{pmatrix}.$$

The subgroup  $$N_1 \varphi = \begin{pmatrix} 1 & 0 \\ U_1 & 1 \end{pmatrix} \leq  \begin{pmatrix} B & 0 \\ T & 1 \end{pmatrix}$$
is in the kernel of this epimorphism, hence  $U_1 \leq \ker \delta$. Therefore, $\delta$ induces a module epimorphism  $\tau :D=T/U_1 \rightarrow D_1$ and a group epimorphism
$$\begin{pmatrix} B & 0 \\ D & 1 \end{pmatrix} \rightarrow \begin{pmatrix} B & 0 \\ D_1 & 1 \end{pmatrix}.$$
By construction  $d(g) \tau =d_1(g), \ g \in G.$ It remains to show that $\tau$ is an isomorphism. From the conditions on the splittings we have two  module epimorphisms  $\sigma :D \rightarrow (B-1) \cdot \MZ B$ and
 $\sigma_1 :D_1 \rightarrow (B-1) \cdot \MZ B$ such that $d(g) \sigma =\overline{g}-1,$  and $d_1(g) \sigma_1 =\overline{g}-1,$ for any $g \in G,$ and such that the following diagram
 \begin{diagram}
D & \rTo^\tau & D_1 \\
\dTo_{\sigma} & & \dTo_{\sigma_1} \\
(B-1) \cdot \MZ B & \rTo^{id} & (B-1) \cdot \MZ B \\
\end{diagram}
is commutative.  Since the restriction of $\tau$ onto $\ker \sigma \cong C$ gives  an isomorphism  $\ker \sigma \rightarrow \ker
\sigma_1$  the epimorphism  $\tau$ is injective,  as claimed.

\end{proof}

The following result easily follows from Lemma \ref{le:splitting}.

\begin{lemma} \label{le:induces-splitting} Let
$$\lambda:G \to \begin{pmatrix} B & 0 \\ D & 1 \end{pmatrix}$$
be a splitting of $G$ over $C$. For a subgroup $A \leq G$ denote by $\overline{A}$ its image in $B$ and by  $D(A)$ the $\MZ \overline{A}$-submodule of $D$ generated by $d(a), \ a \in A.$ Then the restriction $\lambda_A$ of $\lambda$ onto $A$ gives a splitting
 $$\lambda_A: A \to  \begin{pmatrix} \overline{A} & 0 \\ D(A) & 1 \end{pmatrix}$$
 of $A$ over $C \cap A.$
 \end{lemma}

Recall, that if $M$ is a  right module without torsion over a (right) Ore domain  $R$ then the rank $rank_R M$ of $M$ is defined as the cardinality of a maximal system of linearly independent over $R$  elements in $M$. This notion is well-defined, indeed, the ring $R$ embeds into its ring of (right) fractions $K = K(R)$ which is a  division ring. Since the $R$-module $M$ has no $R$-torsion it embeds into its tensor completion $V =  M \bigotimes_{R} K$, which is a vector space over $K$. It is not hard to see that $rank_{R}M = \dim_{K} V$, which is well-defined.

\begin{lemma} \label{le:3} Let  $G$ be a group generated by a finite set of elements $\{ g_1, \ldots ,g_m \}$, $C$ an abelian normal subgroup of a group $G$, and  $B=G/C \neq 1$.  Suppose that the group ring  $\MZ B$ is an Ore domain and $C$, viewed  as  $\MZ B$-module does not have torsion. Then the following hold:
  \begin{itemize}
\item [1)] If
$$\lambda: G \to \begin{pmatrix} B & 0 \\ D & 1 \end{pmatrix}$$
is a splitting of $G$ over $C$, then the module  $D$ has no $\MZ B$-torsion and  ${\rm rank} \ D= {\rm rank} \ C+1 \leq m.$
  \item [2)]  The group $G$ embeds into a wreath product $A_n \wr B$, where $A_n$ is a free abelian group of rank $n = {\rm rank} \ D$.
      \end{itemize}

 \end{lemma}

\begin{proof} To prove 1) notice, first, that the module  $D$ is generated by the set $\{ d(g_1), \ldots ,d(g_m) \}$  (this comes from the homomorphism $\lambda$), so ${\rm rank} \ D \leq m.$  By Lemma \ref{le:splitting} the quotient  $D/C$ is isomorphic to  $(B-1) \cdot ZB$. Clearly,  $(B-1) \cdot ZB$ has rank 1 over $\MZ B$ (as an ideal of $\MZ B$) and it is torsion free since $\MZ B$  does not have zero divisors. It follows that  $D$ also is torsion free (as an extension of a torsion free module $C$ by a torsion free module $D/C$) and ${\rm rank} \ D= {\rm rank} \ C+1$, which proves the statement 1).

To show 2) denote $n = {\rm rank} \ D$.  Since $\MZ B$ is a right (and also left)  Ore domain one can embed $\MZ B$ into its  right ring of  fractions $K$, which is a division ring. In this case the module $D$ embeds into the right vector space $V= D \otimes_{\MZ B}K$ over $K.$ Clearly,  $\dim V=n.$  Let $\{ e_1, \ldots ,e_n \}$ be a basis of  $V.$ One can express the generators $\{
d(g_1), \ldots ,d(g_m) \}$ of  $D$ as linear combinations of the elements of the basis:
$$\begin{matrix} d(g_1)=e_1 \alpha_{11}+ \ldots +e_n \alpha_{1n}, \\
\ldots \\ d(g_m)=e_1 \alpha_{m1}+ \ldots +e_n \alpha_{mn}.
\end{matrix}$$
Notice, that  $K$ is also a left ring of fractions of $\MZ B$ (since $\MZ B$ is also a left Ore domain),  therefore  there exists a non-zero element $u \in \MZ B$ such that  $\alpha_{ij}=u^{-1}v_{ij},$ for suitable $v_{ij} \in \MZ B.$ Clearly, the module  $D$ lies in a free $\MZ B$-module $T$ with basis $\{ t_1=e_1u^{-1}, \ldots ,t_n=e_nu^{-1} \}$. In this event the group $G$ embeds into the group of matrices $$\begin{pmatrix} B & 0 \\ T & 1 \end{pmatrix},$$ which is isomorphic to a wreath product of a free abelian group of rank $n$ and the group $B$.

\end{proof}

\section{Rigid solvable groups}

 The derived (or commutator) subgroup $[G,G]$ of $G$ is denoted by $G^\prime$ or $G^{(1)}$. The $n$-th derived  subgroup $G^{(n)}$ is defined by induction as $G^{(n+1)} = [G^{(n)},G^{(n)}]$, where $G^{(0)} = G$. The group $G$ is solvable of class $c \geq 1$ if $G^{(c)} = 1$, but $G^{(c-1)} \neq 1$.

\begin{definition}
A series of normal subgroups of $G$
\begin{equation} \label{eq:principal}
 G = G_1 > G_2 > \ldots > G_n > G_{n+1} = 1
\end{equation}
is called {\em principal} if the factors $G_i/G_{i+1}$ are abelian groups  which do not have torsion as $\MZ[G/G_i]$-modules, $i = 1, \ldots,n$.
 \end{definition}

\begin{definition}
 Groups with principal series are called {\em rigid}.
 \end{definition}
  It follows immediately from the definition that for a rigid group $G$ with a principal series (\ref{eq:principal}) the group $G/G_i$ also has a principal series
  $$G/G_i > G_2/G_i > \ldots > G_{i-1}/G_i > G_i/G_i = 1.$$
    This allows one to use induction on the length $n$ of the principal series in a rigid group $G$.

\begin{lemma} \label{le:solvability-class}
Let $G$ be a rigid group with a principal series (\ref{eq:principal}). Then
 \begin{itemize}
 \item [1)] $G_i \geqslant G^{(i-1)}$ for $i=1, \ldots ,n+1$;
 \item [2)] $G^{(i-1)} \not \leq G_{i+1}$ for $i=1, \ldots ,n$.
\end{itemize}
In particular, the solvability class of $G$ is equal exactly to $n$.
\end{lemma}

\begin{proof} 1) is the standard fact on solvable groups.

To show 2)  we prove first that $G^{(n-1)} \not \leq G_{n+1} = 1$.  By induction on $n$ (arguing for the group $G/G_n$) we may assume that
$G^{(n-2)} \nleqslant G_n.$ Let $a \in G^{(n-2)} \smallsetminus G_n$ and $b \in G_n, b \neq 1$.  Then  $1 \neq b^{a-1}=b^ab^{-1} = b^{-1}b^a= [b,a] \in G^{(n-2)}.$ Therefore,
 $1 \neq b^{(a-1)^2}=(b^{a-1})^{a-1} = [b,a]^{a-1}= [b,a,a] \in [G^{n-2},G^{n-2}] = G^{(n-1)}$. Hence, $G^{(n-1)} \neq 1$, as claimed.  Now, the inequality $G^{(i-1)} \not \leq G_{i+1}$ for an arbitrary $i = 1, \ldots ,n$ follows from the argument above applied to the group  $G/G_{i+1}$.

\end{proof}

Sometimes we refer to $n$-solvable rigid groups as to $n$-rigid groups.

\begin{lemma} \label{le:definable}
Let $G$ be a rigid group with a principal series (\ref{eq:principal}). Then the following holds for every $i = 1, \ldots, n$:
  \begin{itemize}
 \item [1)] if  $(g_1, \ldots ,g_n)$ is a tuple of elements from $G$ such that  $g_m \in G_m \smallsetminus G_{m+1}$ for every $m = 1, \ldots, n$  then
  $$ G_i= \{ x \in G \ | \ [x,g_i,g_{i+1}, \ldots ,g_n]=1 \}.$$
    \item [2)]
   $$ G_i= \{ x \in G \ | \
[x,G^{(i-1)},G^{(i)}, \ldots ,G^{(n-1)}]=1 \}.$$
\end{itemize}
 \end{lemma}
\begin{proof}
We prove by induction on $i$ that 1) and 2) hold for $G_{n-i}$.  Let $i = 0$. Since $G_n$ is abelian $G_n \subseteq C_G(g_n)$.  To show the converse inclusion assume that $x \in C_G(g_n) \smallsetminus G_n$. Then $g_n^{x-1} = [g,x] = 1$ - contradiction with no $\MZ[G/G_n]$-torsion in $G_n$.  Therefore, $G_n  = C_G(g_n)$ Now by Lemma \ref{le:solvability-class} $G^{(n-1)} \leq G_n$ and $G^{(n-1)} \neq 1$, so $G_n = C_G(G^{(n-1)})$. This proves the base of induction. Suppose now the statements 1) and 2) hold for $i$ and consider the rigid group $G/G_{n-i}$. The argument above shows that $G_{n-i-1} = \{x \in G \mid [x,g_{n-i-1}] \in G_{n-i}\}$. Therefore,
$$x \in G_{n-i-1} \Longleftrightarrow [x,g_{n-i-1}] \in G_{n-i} \Longleftrightarrow [x,g_{n-i-1}, g_{n-i}, \ldots,g_n] = 1,$$
as required in 1). Similarly,
$$x \in G_{n-i-1} \Longleftrightarrow [x,G^{(n-i-2)}] \in G_{n-i} \Longleftrightarrow [x,G^{(n-i-2)}, G^{(n-i-1)}, \ldots,G^{(n-1)}] = 1,$$

\end{proof}

\begin{corollary}
Every rigid group has only one principal series.
\end{corollary}

In the next proposition  we collect  some simple properties of rigid groups.

\begin{proposition}
  The following hold:
   \begin{itemize}
     \item [1)] Rigid groups are torsion-free solvable groups.
     \item [2)] Torsion-free abelian groups are rigid.
     \item [3)] Subgroups of rigid groups are rigid.
     \item [4)] Direct products of two groups, one of which is non-abelian, is not  rigid.
     \item [5)] Non-abelian groups with non-trivial center are not rigid.
   \end{itemize}
\end{proposition}
  \begin{proof}
      To see 1) notice that if (\ref{eq:principal}) is a principal series for a group $G$, then the factors $G_i/G_{i+1}$ are torsion-free as abelian groups. Indeed, if $g\in G_i$ and $g^m \in G_{i+1}$ then $G_i/G_{i+1}$ has $\MZ[G/G_i]$-torsion, since $(gG_{i+1})^{m \cdot 1} = g^mG_{i+1}= G_{i+1}$. Hence $g \in G_{i+1}$. In particular, if $g^m = 1$ then $g = 1$.

      2) If $G$ is abelian then $G > 1$ is a principal series for $G$.

      3) If (\ref{eq:principal}) is a principal series for $G$ and $A \leq G$ then one can construct a normal  series  $A = A\cap G_1 \geq A \cap G_2 \geq \ldots \geq A \cap G_n \geq 1$ and  delete repeating  terms. It is an easy exercise to show that the resulting series is a principal one for $A$.

      4) If $G = A \times B$ with non-abelian $B$, then for any $1 \neq g = (a,b) \in G$ the centralizer $C_G(g)$ contains $(a^\prime,1)$ for some non-trivial $a^\prime \in A$ ($a^\prime  = a$, if $a \neq 1$, and $a^\prime$ is an arbitrary non-trivial from $A$ if $a =1$). In this case $C_G((a^\prime,1))$ contains non-abelian subgroup $(1,B)$. So for any $g \in G$ there is a nontrivial element $h \in C_G(g)$ such that $C_G(h)$ is non-abelian. By Lemma \ref{le:definable} $G$ is non-rigid.

      5) Similar to 4).
  \end{proof}

The following result gives a host of  non-abelian examples  of rigid groups.

\begin{proposition} \label{pr:rigid-examples}
  The following hold:
   \begin{itemize}
     \item [1)] Free solvable groups are rigid.

     \item [2)] If $A$ is torsion-free abelian and $B$ is rigid then the wreath product $A \wr B$ is rigid.
     \end{itemize}
     \end{proposition}
   \begin{proof}  We prove 2) first. Notice, that a torsion-free abelian group $A$ embeds into a divisible torsion-free  abelian group $A_\mathbb{Q}$, here $A_\mathbb{Q}$ is the additive group of a vector space $ A \bigotimes_\mathbb{Z} \mathbb{Q}$  over  rationals $\mathbb{Q}$.  In this case $A \wr B \leq A_\mathbb{Q} \wr B$ and since subgroups of rigid groups are rigid it suffices to prove the result for $A_\mathbb{Q} \wr B$. To this end assume from the beginning that $A$ is a divisible torsion-free abelian group with a maximal linearly independent system of elements $E= \{ e_i \ | \ i \in I \}.$  Then
      $$A \wr B = \begin{pmatrix} B & 0 \\ T & 1
\end{pmatrix},$$
  where $T$ is a right free $QB$-module with basis $E$. If
   $$B=B_1 > B_2 > \ldots B_n >B_{n+1}=1$$
    is a principal series in  $B,$ then
      $$\begin{pmatrix} B_1 & 0 \\ T & 1
\end{pmatrix} > \begin{pmatrix} B_2 & 0 \\ T & 1
\end{pmatrix} > \ldots > \begin{pmatrix} B_n & 0 \\ T & 1
\end{pmatrix} > \begin{pmatrix} 1 & 0 \\ T & 1
\end{pmatrix} >1$$
   is a principal series in $A \wr B$, hence $A \wr B$ is rigid.

Now 1) follows from 2) and the Magnus embedding by induction on the solvability class.

  \end{proof}

For a principal series (\ref{eq:principal}) of a rigid group $G$ we denote the group ring
 $\MZ[G/G_i]$ by $R_i = R_i(G)$ and the right $R_i$-module $G_i/G_{i+1}$ by $T_i = T_i(G)$. Since the group
$G/G_i$ is solvable and torsion-free the group ring $R_i$ is an  Ore domain \cite{KLM}, hence $R_i$ embeds into its ring of fractions $K_i = K_i(G)$ which is a  division ring. Since the $R_i$-module $T_i$ has no $R_i$-torsion it embeds into its tensor completion $V_i = V_i(G) =  T_i \bigotimes_{R_i} K_i$, which is a vector space over $K_i$. Put
$$r_i(G) = \dim_{K_i} V_i = rank_{R_i}T_i, \ \ \ r(G)=(r_1(G), \ldots ,r_n(G)).$$
The tuple $r(G)$ is the {\em principal dimension} of $G$.

\begin{lemma} Let $G$ be an $n$-rigid group. If $G$ is generated by $m$ elements then   $r_1(G) \leq m$ and $r_i(G) \leq m-1 \ (2 \leq i
\leq n)$.
\end{lemma}
  \begin{proof}
  The inequality $r_1(G) \leq m$ is obvious. The inequality $r_n(G) \leq m-1$ follows directly from Lemma \ref{le:3} part 1). The other inequalities follow by induction on $n$ when arguing for the group $G/G_n$.
  \end{proof}

 Now we introduce an iterated wreath product $W(m,n)$ of $n$ free abelian groups $A_m$ of rank $m$. Put $W(m,0) =A_m$ and define $W(m,n)$ by induction
 $W(m,n) = A_m \wr W(m,n-1)$.

 The following result gives a nice characterization of finitely generated rigid groups.

 \begin{proposition} \label{pr:embedding-W(m,n)}
 Let $G$ be an $n$-rigid $m$-generated group. Then $G$ embeds into $W(m,n)$.  Conversely, the group $W(m,n)$ is rigid, so every finitely generated subgroup of $W(m,n)$ is rigid.
 \end{proposition}
 \begin{proof}
  Let $G$ be a rigid $m$-generated group with the principal series (\ref{eq:principal}). By induction on $n$ we may assume that $B = G/G_n$ embeds into $W(m,n-1)$. If in Lemma \ref{le:3}  we put $C = G_n$ then (in the notation of Lemma \ref{le:3}) the rank of the module $D$ is at most $m$, and by the second part of  Lemma \ref{le:3} $G$ embeds into $A_m \wr W(m,n-1) = W(m,n)$.  The converse follows by induction from Proposition \ref{pr:rigid-examples}.
 \end{proof}

 \begin{proposition}
 Finitely generated rigid groups are equationally Noetherian.
 \end{proposition}
   \begin{proof}
   It has been shown in \cite{GR} that the groups $W(m,n)$ are equationally Noetherian. Now the result follows from Proposition \ref{pr:embedding-W(m,n)} and the fact that subgroup of equationally Noetherian groups are equationally Noetherian.
   \end{proof}

In the rest of this section we study behavior of the principal dimension under proper epimorphisms (epimorphisms with non-trivial kernels).
\begin{lemma} \label{le:epimorphism8}
 Let $G$ and $H$ be $n$-rigid groups. If $\varphi :G \rightarrow H$ is a proper epimorphism then  $r(G) > r(H)$ in the left lexicographical order.
\end{lemma}
  \begin{proof}
Let
$$
G=G_1>G_2>\ldots>G_{n+1}=1,
$$
$$
H=H_1>H_2>\ldots>H_{n+1}=1,
$$
be the  principal series of $G$ and $H$. By Lemma \ref{le:solvability-class} the groups $G$ and $H$ are both $n$-solvable. Since $\phi$ is epimorphism it maps $G^{(n-1)}$ onto $H^{(n-1)}$, so  there exists $g\in G^{(n-1)}$ such
that $g^\varphi \in H^{(n-1)} \smallsetminus \{1\}$ (otherwise $H^{(n-1)} = 1$).
  By Lemma \ref{le:definable}  $G_n= C_G(g)$ and $H_n= C_H(g \varphi)$.  Hence $G_n \varphi \leqslant H_n$ and  $\varphi$ gives rise to an epimorphism  $\psi
:G/G_n \rightarrow H/H_n.$  There are two cases to consider.

Case 1. If  $ \psi$ is proper then by induction on $n$ we may assume that $r(G/G_n) > r(H/H_n)$ in the left lexicographical order. Obviously, there are natural isomorphisms  of rings $R_i(G)$ and $R_i(G/G_n)$, as well as the abelian groups $T_i(G)$ and $T_i(G/G_n)$, for $i = 1, \ldots, n-1$, which give rise to an isomorphism of the modules  $T_i(G)$ and $T_i(G/G_n)$ (under the proper identification of scalars). This implies that $r_i(G) = r_i(G/G_n), r_i(H) = r_i(H/H_n)$ for $i = 1, \ldots, n-1$, hence $r(G) > r(H)$.

Case 2. If  $\psi$ is an isomorphism then $r(G/G_n) =  r(H/H_n)$ and  $\ker \varphi \leqslant G_n$. In this case, as was argued above, $r_i(G) = r_i(H)$ for $i = 1, \ldots, n-1$. Moreover,  $H_n$, viewed as  $Z[G/G_n]$-module (upon identification of coefficients via $\psi$), is isomorphic to the factor module $G_n/ \ker \varphi.$ Therefore, $r_n(G) > r_n(H),$ so $r(G) > r(H)$, as claimed.
\end{proof}

\begin{lemma} Let $A$ be an $n$-rigid subgroup of an $n$-rigid group $G$. If $\varphi :G \rightarrow A$ is an $A$-epimorphism then $G_i \varphi = A \cap G_i$  for every $i = 1, \ldots, n$.
 \end{lemma}
 \begin{proof} Put $A_i = A \cap G_i$. We show first that  $G_n \varphi = A_n$.
Let  $1 \neq a \in A_n.$ By Lemma \ref{le:definable} $G_n=C_G(a).$ Obviously, $C_G(a)
\varphi \leqslant C_A(a)=A_n.$ On the other hand $A_n \varphi =A_n$, so  $G_n \varphi =A_n$, as claimed. Now $\varphi$ gives rise to an epimorphism $G/G_n \rightarrow A/A_n$ and induction finishes the proof.
\end{proof}

\section{Discrimination and rigid groups}
  \label{se:disc}

In this section  we study properties that are inherited under discrimination by rigid groups
 with an eye on applications to the coordinate groups of irreducible varieties.

\begin{lemma} \label{le:disc10} Let $A$ be an $n$-rigid group with the principal series
\begin{equation} A=A_1 > A_2 > \ldots > A_n > A_{n+1}=1
\end{equation}
and $G$  a group, containing $A$ as a subgroup. If  $\Phi$ is a set of $A$-epimorphisms from $G$ into $A$  that $A$-discriminates $G$ then:
  \begin{itemize}
   \item [1)] $G$ is an $n$-rigid group;
    \item [2)] if $G = G_1 > G_2 > \ldots > G_n > G_{n+1}=1 $ is the principle series of $G$ then $G_i \varphi = A_i$ for every $\varphi \in \Phi$, $ i = 1, \ldots,n+1$;
        \item [3)]   the group $A/A_i$ discriminates the group $G/G_i$, $ i = 1, \ldots,n+1$,  by the set of epimorphisms induced from $\Phi$.
\end{itemize}
 \end{lemma}
\begin{proof}
Let  $a \in A_n, a \neq 1.$ By Lemma \ref{le:definable} $A_n = C_A(a)$. Put  $G_n =C_G(a).$ We claim that $G_n$ is a normal abelian subgroup of $G$. Indeed,  observe, first, that $G_n \varphi = A_n$ for every $\varphi \in \Phi$ (since $\varphi(a) = a$), so the set of the restrictions on $G_n$ of epimorphisms from $\Phi$ $A_n$-discriminates $G_n$ into $A_n$, in particular, $G_n$ is abelian. Secondly, $G_n$ is a normal subgroup of $G$. To see this let  $h \in G_n, \ g \in G$ and  $h^g \notin G_n,$ so  $[h^g,a] \neq 1.$ Then there is $\varphi \in \Phi$ such that
$[h^g,a] \varphi \neq 1$, which implies  $[(h \varphi)^{g \varphi},a] \neq 1$ in $A$ - contradiction with normality of $A_n$ in $A$. This proves the claim.

Notice now that every epimorphism  $\varphi \in \Phi$ induces an epimorphisms $\bar{\varphi}$ of quotient groups $G/G_n \to A/A_n$ that is identical on $A/A_n$ (which is viewed as a subgroup of $G/G_n$ via the natural embedding). Denote the set of such epimorphisms by $\bar{\Phi}$. We claim that $\bar{\Phi}$ $A/A_n$-discriminates $G/G_n$ into $A/A_n$.
Indeed, let  $g_1G_n, \ldots ,g_mG_n$ be pairwise distinct elements from $G/G_n$. Then  $[g_ig_j^{-1},a] \neq 1$ for  $i \neq j.$ Since $\Phi$ discriminates $G$ into $A$  there exists an epimorphism  $\varphi \in \Phi$ such that
$[g_ig_j^{-1},a] \varphi \neq 1.$ It follows that $g_1 \varphi \cdot A_n, \ldots ,g_m \varphi \cdot A_n$ are pairwise distinct elements of $A/A_n$,  so $\bar{\varphi} \in \bar{\Phi}$ discriminates the elements $g_1G_n, \ldots ,g_mG_n$ into $A/A_n$, as claimed.

Now by induction on $n$ we may assume that the statements of the lemma hold for the group  $G/G_n$.  It remains to be seen that $G_n$, as an  $\MZ [G/G_n]$-module does not have torsion. Let  $h \in G_n, h \neq 1,$ and  $u \in \MZ [G/G_n], u \neq 0$.  We need to show that $h^u \neq 1.$ Let  $u= \alpha_1 \cdot g_1G_n+ \ldots + \alpha_m \cdot g_mG_n,$ where  $ \alpha_i \in \MZ, \alpha_i \neq 0$, and $g_1G_n, \ldots ,g_mG_n$ are pairwise distinct elements of $G/G_n.$ As we have seen above there exists  $\varphi \in \Phi$ such that $h \varphi \neq 1$ and  $g_1 \varphi \cdot A_n, \ldots ,g_m \varphi \cdot A_n$ are pairwise distinct elements of $A/A_n.$ Since  $A_n$, as an  $\MZ[A/A_n]$-module does not have torsion then  $(h \varphi)^{\alpha_1 \cdot g_1 \varphi \cdot A_n+
\ldots + \alpha_m \cdot g_m \varphi \cdot A_n} \neq 1.$ Hence $h^u \neq 1$, as required.
This proves the lemma.
\end{proof}

\begin{lemma} \label{le:disc11} Let $A$ be an $n$-rigid group
and $G$  a group, containing $A$ as a subgroup. If  $\Phi$ is a set of $A$-epimorphisms from $G$ into $A$  that $A$-discriminates $G$ then $G$ is $n$-rigid group and the following conditions hold:
  \begin{itemize}
   \item [1)] every system of elements from $T_i(A)$ linearly independent over $\MZ [A/A_i]$ is also linearly independent over  $\MZ [G/G_i]$.  In particular, $r_i(G) \geq r_i(A),$ $ i = 1, \ldots,n$;
    \item [2)] if  $G$ is $m$-generated, as an $A$-group, and $r_i(A)$ is finite  then $r_i(G)-r_i(A) \leq m$, $ i = 1, \ldots,n$.
\end{itemize}
 \end{lemma}
\begin{proof}
To prove the lemma it suffices to show that the statements 1) and 2) hold for $i = n$,   the rest follows from Lemma \ref{le:disc10} by induction on $n$ for the groups $G/G_n$ and $A/A_n$.

To see that  $r_n(G) \geq r_n(A)$ choose some  system of elements $\{t_1, \ldots ,t_s \} \in A_n$ which is linearly independent over the ring $\MZ [A/A_n]$. It suffices to show that this system is also linearly independent over $\MZ [G/G_n]$. Suppose to the contrary that  $t_1u_1+ \ldots +t_su_s=0,$ where $u_i \in \MZ[G/G_n]$ and, say,  $u_1 \neq 0.$ As in Lemma \ref{le:disc10} there exists  $\varphi \in \Phi,$ such that  $u_1\bar{\varphi} \neq 0$ in the ring  $\MZ[A/A_n]$. It follows that  $t_1(u_1 \bar{\varphi})+ \ldots +t_s(u_s \bar{\varphi})=0$  - contradiction with linear independence of $\{ t_1, \ldots ,t_s \}$ over  $\MZ[A/A_n]$. This proves 1).

To prove that $r_n(G)-r_n(A) \leq m$ consider a set of generators $\{ g_1, \ldots ,g_m \}$ of $G$, as an $A$-group.   Let
$$\lambda: G \to \begin{pmatrix} B & 0 \\ D
& 1 \end{pmatrix}$$
be a splitting of $G$ over $G_n$, where  $B=G/G_n.$ In the notation of Lemma \ref{le:splitting}  $\MZ B$-module $D$
is generated by elements $\{ d(g_1), \ldots ,d(g_m),d(a) \ | \ a \in A \}.$ Denote by $D(A)$ the $\MZ[A/A_n]$-submodule of $D$ generated by the elements  $\{ d(a) \ | \ a \in A \}.$ By Lemma \ref{le:induces-splitting} the group
$$\begin{pmatrix} A/A_n & 0 \\ D(A) & 1 \end{pmatrix}$$
gives a splitting of $A$ over  $A_n.$ By Lemma \ref{le:3} ${\rm rank} \ D(A)=r_n(A)+1.$ Let now  $r_n(A)=s$ and  $\{ h_1, \ldots ,h_s \}$ a maximal linearly independent over $\MZ[A/A_n]$   system of elements from $A_n$.  Take an arbitrary  $h_0 \in D(A) \setminus A_n.$ The argument above show that the elements  $\{ h_1, \ldots ,h_s \}$ are linearly independent over $\MZ B.$  Moreover, we claim that the system of elements  $\{ h_0, h_1, \ldots ,h_s \}$ is also linearly independent over  $\MZ B.$ Indeed, otherwise, there exists a linear combination
$h_0u_0+ h_1u_1+ \ldots +h_su_s=0,$ where  $u_i \in \MZ B, \ u_0 \neq 0.$ Applying to the both sides the epimorphism
$\sigma :D \rightarrow (B-1)\MZ B$ one gets $h_0 \sigma \cdot u_0 =0$  (since $h_i \in A_n, i \neq 0$ and $B = G/G_n$) -contradiction with  $h_0 \sigma \neq 0,\ u_0 \neq 0.$

Now consider the  $\MZ B$-module  $D(A) \cdot ZB,$ generated by the set of elements  $D(A).$ Since  $\{ h_0, h_1, \ldots ,h_s \}$ is a maximal linearly independent system of elements of
 $\MZ[A/A_n]$-module $D(A)$ , which is also linearly independent over  $\MZ B,$ then ${\rm
rank} \ (D(A) \cdot ZB)=s+1.$ Observe, that
$$D=d(g_1) \cdot ZB+ \ldots +d(g_m) \cdot ZB+ D(A) \cdot ZB,$$
hence  ${\rm rank} \ D \leq m+s+1= m+r_n(A)+1.$ It follows that  $r_n(G)= {\rm rank} \ D -1 \leq m +r_n(A_n)$, as required.
\end{proof}

A slight modification of the proof of the statement 2) above gives the following result.

\begin{remark}
 An analog of the statement 2) in Lemma \ref{le:disc11} also holds when the ranks $r_i(A)$ are allowed to be infinite. In this case we claim that the rank $T_i(G)$ {\em modulo} $T_i(A)$ is at most $m$, or equivalently, the rank of the $\MZ[G/G_i]$-module $T_i(G)/T_i(A)\MZ[G/G_i]$ is at most $m$, where  $T_i(A)\MZ[G/G_i]$ is the $\MZ[G/G_i]$-submodule of $T_i(G)$ generated by $T_i(A)$.
\end{remark}

\section{Main Theorems}
In this section we prove the main theorems.

\begin{theorem} \label{th:1} Let $A$ be a finitely generated  $n$-rigid group, $S$ a non-empty irreducible algebraic set from $A^m$, $G=\Gamma (S)$ the coordinate group of $S$. Then $G$ is an $n$-rigid group and $0 \leq r_i(G)-r_i(A) \leq m$ for all $i=1, \ldots n$.
 \end{theorem}
 \begin{proof}
   By Theorem \ref{coord-groups} the coordinate group $G$ is an $m$-generated $A$-group which  is $A$-discriminated by $A$. By Lemma \ref{le:disc10} $G$ is also $n$-rigid group and by Lemma \ref{le:disc11} $r_i(G)-r_i(A) \geq 0.$ Since  $G$, as an  $A$-group, is $m$-generated then  by Lemma \ref{le:disc11}
$r_i(G)-r_i(A) \leq m,$ which proves the theorem.
 \end{proof}

In the notation of Theorem \ref{th:1} define  $\alpha (S)=r(G)-r(A).$

\begin{theorem} \label{th:2} Let  $A$ be a finitely generated  $n$-rigid group. Then
   \begin{itemize}
   \item [1)] if  $P \subsetneqq S $  are irreducible algebraic sets from $A^m$ then  $\alpha (S) > \alpha (P).$
    \item [2)] The Krull dimension of any irreducible algebraic set from $A^m$ does not exceed  the number $(m+1)^n$. In particular, the Krull dimension of $A^m$ is finite for every $m$.
        \end{itemize}
  \end{theorem}

\begin{proof} Let in the notation above  $G= \Gamma (S)$  and $H= \Gamma (P).$ By Theorem \ref{th:1} $G$ ш $H$ are $n$-rigid  $A$-groups. Since  $P \subseteq S$, then (see Section \ref{se:Krull}) there exists a proper  $A$-epimorphism $G \rightarrow H$. By Lemma \ref{le:epimorphism8}  $r(G)>r(H)$  in the left lexicographical order. It follows that $\alpha (S)=r(G)-r(A)> \alpha
(P)=r(H)-r(A)$, which proves 1).

To see 2) observe from Theorem \ref{th:1}  that  for any irreducible algebraic set $S \subseteq A^m$ one has $\alpha(S) \leq (m, \ldots ,m)$, where the tuple has the length $n$.   Hence there are at most $(m+1)^n$ such tuples $\alpha(S)$, so the Krull dimension of any irreducible algebraic set in $A^m$ is bounded from above by $(m+1)^n$, as claimed.

\end{proof}

\end{document}